\documentclass[12pt]{article}
\usepackage{amsmath}
\usepackage{amssymb}

\usepackage{mathtools}
\usepackage{amsthm}
\usepackage[usenames]{color}

\usepackage{graphicx}

\newcommand{\1}[1]{{\mathbf 1}{\{#1\}}}

\newcommand{\vr}{\varrho}

\newcommand{\Z}{{\mathbb Z}}

\newcommand{\HH}{{\mathcal H}}
\newcommand{\R}{{\mathbb R}}

\newcommand{\hg}{\hat{g}}

\newcommand{\s}{{\widehat S}}
\newcommand{\hN}{{\widehat N}}

\newcommand{\hG}{\widehat{G}}

\newcommand{\hl}{{\hat\ell}}

\let\phi=\varphi

\newcommand{\E}{{\mathbb E}}

\newcommand{\diam}{{\mathop{\mathrm{diam}}}}

\newcommand{\hM}{\widehat{M}}

\newcommand{\IP}{{\mathbb P}}
\newcommand{\IE}{{\mathbb E}}

\newcommand{\hP}{\widehat{P}}

\newcommand{\capa}{\mathop{\mathrm{cap}}}

\newcommand{\hcapa}{\mathop{\widehat{\mathrm{cap}}}}
\newcommand{\hm}{\mathop{\mathrm{hm}}\nolimits}
\newcommand{\hhm}{\mathop{\widehat{\mathrm{hm}}}\nolimits}

\newcommand{\hEs}{\mathop{\widehat{\mathrm{Es}}}\nolimits}

\newcommand{\dist}{\mathop{\mathrm{dist}}}

\newcommand{\htau}{\widehat{\tau}}

\newcommand{\B}{{\mathsf B}}

\allowdisplaybreaks

\newtheorem{theo}{Theorem}[section]
\newtheorem{lem}[theo]{Lemma}

\newtheorem{prop}[theo]{Proposition}

\title{Conditioned two-dimensional simple random walk: 
Green's function and harmonic measure}
\author{Serguei~Popov}

\begin{document}

\maketitle

{\footnotesize 
\noindent Department of Statistics, Institute of Mathematics,
 Statistics and Scientific Computation, University of Campinas --
UNICAMP, rua S\'ergio Buarque de Holanda 651,
13083--859, Campinas SP, Brazil\\
\noindent e-mail: \texttt{popov@ime.unicamp.br}

}

\begin{abstract}
We study the Doob's $h$-transform of the two-dimensional simple random walk
with respect to its potential kernel, which can be thought of as
 the two-dimensional simple random walk
conditioned on never hitting the origin. 
We derive an explicit formula for the Green's function
of this random walk, and also prove a quantitative 
result on the speed of convergence of the (conditional)
entrance measure to the harmonic measure
(for the conditioned walk) on a finite set.
\\[.3cm]\textbf{Keywords:} transience, 
 Doob's $h$-transform, entrance measure
\\[.3cm]\textbf{AMS 2010 subject classifications:}
Primary 60J10. Secondary 60G50, 82C41.
 
\end{abstract}

\section{Introduction and main results}
\label{s_introres}
In this paper, we derive some basic ``potential-theoretic''
results for the two-dimensional simple random walk
conditioned on never hitting the origin, namely
\begin{itemize}
 \item[(i)] we obtain an explicit expression for its
 Green's function;
 \item[(ii)] we prove a \emph{quantitative} result on the convergence
of the conditional entrance measure to a finite set~$A$ starting from
a distant site to the harmonic measure~$\hhm_A$
of the conditioned walk.
\end{itemize}
As usual, results on the Green's function and the control
of the entrance measure are important tools in the study
of properties of trajectories of the random walk.
Regarding the above item~(i), one can note that the expression
for that Green's function was \emph{implicitly} present in 
the paper~\cite{ComPopVac16}; however, here we give a more
clean-and-straightforward derivation. 
Regarding item~(ii), it is important to keep in mind
that the harmonic measure~$\hhm_A$ with respect to 
the conditioned walk is generally different from the 
harmonic measure~$\hm_A$ with respect to 
the simple random walk.

The two-dimensional simple random walk conditioned on never 
hitting the origin is the main ingredient in the construction
of the two-dimensional random interlacements introduced in~\cite{ComPopVac16}
and further studied in~\cite{ComPop16,Rod19} (by
its turn, it is an 
extention of classical random intelacement 
model~\cite{CerTei12,DreRatSap14,Szn10}
to two dimensions).
It then became evident that 
the conditioned walk (denoted by~$\s$ in this paper) is an interesting
object on its own. 
We list its basic properties later in this section, but
let us mention a few more advanced but \emph{surprizing}
ones~\cite{GanPopVac19,PopRolUng19}:
\begin{itemize}
 \item Although the conditioned walk itself is transient,
it hits \emph{any} fixed infinite subset of~$\Z^2$ infinitely many
times a.s.\ (i.e., any infinite set is recurrent for~$\s$).
 \item Fix a \emph{typical} large set (e.g., a large disk or rectangle,
 a long segment, etc.). Then, the proportion of the sites of 
 this set which are eventually visited by~$\s$ is close
 in distribution to the Uniform law on the unit interval.
 \item Let $M_n := \min_{m \geq n} \|\s_m\|$ be the 
 future minimal distance to the origin for the conditioned
 walk (here and in the sequel $\|\cdot\|$ stands for the 
 Euclidean norm). 
 By transience, it is clear that $M_n\to \infty$ a.s.,
but this divergence occurs in a ``highly irregular'' way:
for any~$\delta>0$, it holds that $M_n \leq n^{\delta }$
and $M_k \ge k^{\frac{1}{2}-\delta}$, for infinitely 
many $n,k$.
 \item despite transience, two independent copies of the 
 $\s$-walks (starting from sites of the same parity) will meet
 infinitely many times a.s.; the same will happen for a pair
 of independent $\s$-walk and simple random walk.
\end{itemize}

Also, the content of this paper is intended to be a 
part of the book~\cite{2srw}, that the author is currently 
working on\footnote{comments on~\cite{2srw} will be, of course,
very welcome}.

Now we give formal definitions.
We will repeatedly use the big-O notation: $f(x)=O(g(x))$
as $x\to a$ means that 
$\limsup_{x\to a}\big|\frac{f(x)}{g(x)}\big|<\infty$,
where $a\in\R\cup\{\infty\}$.
In the following, we denote by $(S_n,n\geq 0)$ the simple random
walk in~$\Z^2$.
Let us define its \emph{potential kernel} $a(\cdot)$ by
\begin{equation}
\label{def_a(x)}
a(x) = \sum_{k=0}^\infty\big(\IP_0[S_k=0]-\IP_x[S_k=0]\big),
 \qquad x\in \Z^2.
\end{equation}
By definition, it holds that $a(0)=0$, and one can show that
the above series converges and that the resulting value
is strictly positive for all~$x\neq 0$ (here and in the sequel
we refer to Section~4.4 of~\cite{LawLim10}). 
Also, the function~$a$ is harmonic outside the origin,
i.e.,
\begin{equation}
\label{a_harm}
  a(x) = \frac{1}{4}\sum_{y\sim x}a(y)\quad \text{ for all }
 x\neq 0.
\end{equation}
It is possible to prove that, 
as $x\to\infty$, 
\begin{equation}
\label{formula_for_a}
 a(x) = \frac{2}{\pi}\ln \|x\| + \gamma' + O(\|x\|^{-2}) ,
\end{equation}
where,
 $\gamma'=\pi^{-1}(2\gamma+\ln 8)$  
with $\gamma=0.5772156\dots$ the Euler-Mascheroni 
constant,
cf.\ Theorem~4.4.4 of~\cite{LawLim10}.

Let $\tau$ and $\tau^+$ be the entrance and the hitting
times for the simple random walk:
\begin{align*}
 \tau_A  &= \min\{n \geq 0: S_n\in A\},
 \\ \intertext{ and } 
 \tau^+_A   &= \min\{ n \geq 1: S_n \in A\},
\end{align*}
and we let $\tau_x:=\tau_{\{x\}}$, $\tau^+_x:=\tau^+_{\{x\}}$.
We also use the notations $\htau_A,\htau^+_A,\htau_x,\htau^+_x$
when~$S$ is substituted by~$\s$ in the above definitions.
Observe that the harmonicity of~$a$ outside the origin
immediately implies 
that the process $a(S_{k\wedge \tau_0 })$ is a martingale.
We will repeatedly use this fact in the sequel.
What we will also repeatedly use, is that, due to~\eqref{formula_for_a},
\begin{equation}
\label{diff_a_logs}
 a(x+y)-a(x) = O\big(\tfrac{\|y\|}{\|x\|}\big)
\end{equation}
for all $x,y\in\Z^2$ such that (say) $\|x\|>2\|y\|$.

With some (slight) abuse of notation, we also consider
the function 
\[
a(r)=\frac{2}{\pi}\ln r + \gamma' 
\]
of a \emph{real} argument~$r\geq 1$. Note that, in general, $a(x)$
need not be equal to $a(\|x\|)$, although they are 
of course quite close for large~$x$. 
The advantage of using this notation is e.g.\ that, 
due to~\eqref{formula_for_a} and~\eqref{diff_a_logs}, we may write
(for fixed~$x$ or at least~$x$ such that $2\|x\|\leq r$)
\begin{equation}
\label{real_a}
 \sum_{y\in\partial \B(x,r)} \nu(y)a(y) = a(r) 
+ O\big(\tfrac{\|x\|\vee 1}{r}\big)
 \qquad \text{ as $r\to\infty$ }
\end{equation}
for \emph{any} probability measure~$\nu$ on $\partial \B(x,r)$.

Now we define the main object of study in this paper,
the Doob's $h$-transform of the two-dimensional simple
random walk with respect to the potential kernel~$a$,
which can be informally seen as ``the simple random walk
conditioned on never hitting the origin''.
It is a Markov chain $(\s_n, n\geq 0)$
on~$\Z^2\setminus \{0\}$, and
with transition probability matrix  
\begin{equation}
\label{def_hatS}
 \hP(x,y) = \begin{cases}
           \displaystyle\frac{a(y)}{4a(x)}, & \text{ if }
            x\sim y, x\neq 0,\\
             0,\phantom{\int\limits^A} & \text{ otherwise}
          \end{cases}
\end{equation}
(note that~\eqref{a_harm} implies that the above is 
a stochastic matrix indeed).

Following~\cite{ComPopVac16},
we summarize the basic properties of the random walk~$\s$:
\begin{itemize}
 \item[(i)] The walk~$\s$ is reversible, with the reversible
 measure~$\mu(x)=a^2(x)$.
 \item[(ii)] In fact, it can be represented as a random walk
 on the two-dimensional lattice with the set of conductances
 $\big(a(x)a(y), x,y\in \Z^2, x\sim y\big)$.
 \item[(iii)] 
The process
 $1/a(\s_{n\wedge \htau_{\mathcal{N}}})$ is a 
martingale, where~$\mathcal{N}$ is the set of the four
 neighbours of the origin.
 \item[(iv)]
As a consequence (see e.g.\ Theorem~2.5.8 of~\cite{MenPopWad17}),
the walk $\s$ is transient.
\item[(v)] 
Moreover, for all $x\neq 0$ 
\begin{equation} 
\label{escape_from_site}
\IP_x\big[\htau^+_x<\infty\big] = 1-\frac{1}{2a(x)},
\end{equation}
and for all $x\neq y$, $x,y\neq 0$ 
\begin{equation}
\label{not_hit_site}
\IP_x\big[\htau_y<\infty\big] = 
\IP_x\big[\htau^+_y<\infty\big] = \frac{a(x)+a(y)-a(x-y)}{2a(x)}.
\end{equation}
\end{itemize}

Next, we define the Green's function of the conditioned walk~$\s$
 in the usual way: for $x,y\in \Z^2\setminus\{0\}$
\begin{equation}
\label{df_GF_conditional} 
\hG(x,y) = \IE_x \sum_{k=0}^{\infty} \1{\s_k=y},
\end{equation}
i.e., $\hG(x,y)$ is the mean number of visits to~$y$ starting
from~$x$ (counting the possible ``visit'' at time~$0$).
We are able to \emph{calculate}
this function in terms of the potential kernel~$a$:
\begin{theo}
\label{t_GF_conditional}
For all $x,y\in \Z^2\setminus\{0\}$ it holds that 
\begin{equation}
\label{eq_GF_conditional} 
\hG(x,y) = \frac{a(y)}{a(x)}\big(a(x)+a(y)-a(x-y)\big).
\end{equation}
\end{theo} 

It is easy to observe that~\eqref{eq_GF_conditional}
actually follows from~\eqref{escape_from_site}--\eqref{not_hit_site}
(which appear already in~\cite{ComPopVac16}).
Indeed, since the total number of visits to~$y$ has
Geometric distribution with success parameter~$\IP_y[\htau^+_y=\infty]$
under~$\IP_y$, we have
\[
 \hG(x,y) = \IP_x[\htau_y<\infty] \times \frac{1}{\IP_y[\htau^+_y=\infty]},
\]
which indeed leads to~\eqref{eq_GF_conditional}.
However, the proof of~\eqref{escape_from_site}--\eqref{not_hit_site}
in~\cite{ComPopVac16} is somewhat involved and not very intuitive.
Here, we take the ``classical'' route of first obtaining
the expression for the Green's function; then, it is straightforward
to derive~\eqref{escape_from_site}--\eqref{not_hit_site} 
from it in the usual way.

Next, we define the capacity and the harmonic measure
for the conditioned walk.
As in~\cite{ComPopVac16} (formulas (13)--(14)),
the capacity~$\hcapa(\cdot)$ with respect to the conditioned walk
is defined by
\begin{equation}
\label{df_hcapa}
 \hcapa(A) = \sum_{y\in A} a^2(y)\hEs_A(y),
\end{equation}
where
\begin{equation}
\label{def_hEs_A}
  \hEs_A(x)=\IP_x[\htau^+_A=\infty]\1{x\in A}
\end{equation}
is the escape probability from~$x\in A$ (again, with respect
to the conditioned walk).
Being $\capa(\cdot)$ the capacity for the two-dimensional
simple random walk (see Section~6.6 of~\cite{LawLim10}),
it holds that
\begin{equation}
\label{eq_capa=hcapa}
\hcapa(A)=\capa\big(A\cup\{0\}\big)
\end{equation}
for all $A\subset\Z^2\setminus\{0\}$, see 
Proposition~2.2 of~\cite{ComPopVac16}.
Then, we also have the usual
(as e.g.\ in Proposition 4.6.4 of~\cite{LawLim10})
transient-case relation
\begin{equation}
\label{Green_hEs}
 \IP_x[\htau_A<\infty] = \sum_{y\in A} \hG(x,y) \hEs_A(y)
   = \sum_{y\in A} \frac{\hG(x,y)}{a^2(y)} \times a^2(y) \hEs_A(y).
\end{equation}

Now, let us define the harmonic measure with respect
to the conditioned walk in the following way:
\begin{equation}
\label{df_hhm_A}
 \hhm_A(y) = \frac{a^2(y)\hEs_A(y)}{\hcapa(A)} .
\end{equation}
(One can also note the following fact about the relation
of~$\hhm_A$ to the harmonic measure~$\hm_A$ with respect
to the simple random walk:
$\hhm_A$ is~$\hm_A$ biased by~$a$ by~(13) and~(15) of~\cite{ComPopVac16}; 
 however, we will not need this fact in this paper.)
%
%
For $x\in \Z^2$ and $A\subset \Z^2$ let us define
$\dist(x,A)=\inf_{y\in A}\|x-y\|$;
 also, denote $\diam(A)=\sup_{x,y\in A}\|x-y\|$.
We then argue that $\hhm_A$ is indeed ``the conditional
entrance measure from infinity'', and, moreover,
we obtain a quantitative result for the difference between the entrance
measure of the conditioned walk to a finite set and the 
harmonic measure~$\hhm_A$ on this set:
\begin{theo}
\label{t_conv_harm_hatS}
Assume that $A\subset \Z^2\setminus\{0\}$ is finite 
and $x\notin A$ is such that 
$\dist(x,A)\geq 12(\diam(A)+1)$. 
For all~$y\in A$, we have
\begin{equation}
\label{eq_conv_harm_hat}
\IP_x[\s_{\htau_A}=y\mid \htau_A<\infty]
= \hhm_A(y)
\big(1+ O\big(\tfrac{\diam(A)}{\dist(x,A)}\big)\big).
\end{equation}
\end{theo}
Let us also mention that the term $O\big(\tfrac{\diam(A)}{\dist(x,A)}\big)$
is the same that one obtains in the corresponding result
for the simple random walk, see Theorem~3.17 of~\cite{2srw}.

\section{Some auxiliary definitions and results}
\label{s_aux}
We use the following notations
for $A\subset \Z^2$: $A^c=\Z^2\setminus A$ is the complement 
of~$A$, $\partial A=\{x\in A : \text{ there exist }y\in A^c
\text{ such that }x\sim y\}$ is the boundary of~$A$,
and $\partial_e A = \partial A^c$ is the external boundary 
of~$A$; also, $\B(x,r)=\{y:\|y-x\|\leq r\}$ is the ball (disk)
in 
$\Z^2$ and $\B(r)$ stands for $\B(0,r)$.
Let us recall a few facts about hitting probabilities
of (conditioned or not) simple random walks, and about
the relationship between~$S$ and~$\s$.
(In the three following results, the error term $O(\cdot)$ is uniform 
with respect to the starting point~$x$.)
\begin{lem}
\label{l_escape_origin_d=2}
Let $y\in\Z^2$ and $r\geq \|y\|$ (so that $0\in \B(y,r)$). 
Assume that $x\in \B(y,r)$ and $x\neq 0$.
Then
\begin{equation}
 \label{nothit_0_dim2}
\IP_x\big[\tau_{\partial\B(y,r)}<\tau^+_0\big] = 
\frac{a(x)}{
 a(r)+ O\big(\frac{\|y\|+1}{r}\big)},
\end{equation}
as $r\to \infty$.
The above also 
holds with $\tau_{\partial_e\B(y,r)}$ on the place
of $\tau_{\partial \B(y,r)}$. 
\end{lem}
\begin{proof}
 This is (part of) Lemma~3.1 of~\cite{ComPopVac16}.
\end{proof}

\begin{lem}
\label{l_escape_hatS}
Assume $r \geq 1$ and $\|x\| \geq r+1$. We have
\begin{equation}
\label{eq_escape_hatS}
\IP_x\big[\htau_{\B(r)}=\infty\big] = 
1-\frac{a(r)+O(r^{-1})}{a(x)}.
\end{equation}
\end{lem}
\begin{proof}
This is Lemma~3.4 of~\cite{ComPopVac16}.
\end{proof}

For $D \subset \Z^2$, 
let~$\Gamma^{(x)}_D$ be the set of all 
nearest-neighbour finite
trajectories that start at~$x\in D\setminus\{0\}$ 
and end when entering~$\partial D$ for the first time;
 denote also~$\Gamma^{(x)}_{y,R}=\Gamma^{(x)}_{\B(y,R)}$.
 For~$\HH\subset \Gamma^{(x)}_D$ write
 $S\in \HH$ (respectively, $\s\in \HH$) if there exists~$k$ such that 
$(S_0,\ldots,S_k)\in \HH$ (respectively, $(\s_0,\ldots,\s_k)\in \HH$). 
In the next result we show that 
$\IP_x\big[S\in \cdot\mid\tau_0>\tau_{\partial \B(R)}\big]$ 
and $\IP_x[\s\in \cdot\,]$ are 
almost indistinguishable on~$\Gamma^{(x)}_{0,R}$ (that is, the conditional
law of~$S$ almost coincides with the unconditional law 
of~$\s$), which justifies the intuition that~$\s$ is~$S$
conditioned on never hitting the origin.
\begin{lem}
\label{l_relation_S_hatS}
Let $x\in \B(R)\setminus\{0\}$, and 
assume $\HH\subset \Gamma^{(x)}_{0,R}$.
We have 
\begin{equation}
\label{eq_relation_S_hatS}
\IP_x\big[S\in \HH\mid \tau_0>\tau_{\partial \B(R)}\big]
 =\IP_x\big[\s \in \HH\big] \big(1+O((R \ln R)^{-1})\big).
\end{equation}
\end{lem}
\begin{proof}
This is Lemma~3.3(i) of~\cite{ComPopVac16}.
\end{proof}

Now, we work with the Green's function of the conditioned walk.
Observe that, by~\eqref{eq_GF_conditional}, we have
\[
 \frac{\hG(x,y)}{a^2(y)} = \frac{\hG(y,x)}{a^2(x)}
  = \frac{a(x)+a(y)-a(x-y)}{a(x)a(y)},
\]
and so it is natural to introduce new notation
$\hg(x,y) = \frac{\hG(x,y)}{a^2(y)} = \hg(y,x)$
for the ``symmetrized'' conditional Green's function.

At this point, let us recall that the function $1/a(\cdot)$ is
 harmonic\footnote{with respect to the conditioned walk}
 on $\Z^2\setminus(\mathcal{N}\cup\{0\})$, and
 observe that  
the Green's function $\hG(\cdot,y)$ 
is harmonic on $\Z^2\setminus\{0,y\}$
(this is an immediate consequence of the total expectation formula).
It turns out that this ``small'' difference will be quite important:
indeed, the latter fact will be operational in some places below, 
for applying the Optional Stopping Theorem in some particular settings.
For future reference, we formulate the above fact in the equivalent
form:
\begin{prop}
\label{p_hG_mart}
For any $y\in\Z^2\setminus\{0\}$
it holds that the process 
 $(\hG(\s_{n\wedge\htau_y},y), n\geq 0)$ is a martingale,
as well as the process $(\hg(\s_{n\wedge\htau_y},y), n\geq 0)$.
 Moreover, let us define 
\begin{equation}
\label{def_hl}
\hl(x,y) = 1 + \frac{a(y)-a(x-y)}{a(x)}
 = \frac{\hG(x,y)}{a(y)} .
\end{equation}
Then the process 
$(\hl(\s_{n\wedge\htau_y},y), n\geq 0)$ is a martingale.
\end{prop}
By the way, notice that
\begin{equation}
\label{hl_to_0}
 \lim_{\|x\|\to \infty} \hl(x,y) = 0
\end{equation}
for any fixed~$y$, so the last process is a
``martingale vanishing at infinity'',
which makes it more convenient for applications
via the Optional Stopping Theorem
(so this is why we kept ``$1+$'' in~\eqref{def_hl}).

Next, we need a few technical estimates on the function~$\hg$.
\begin{lem}
\label{l_order_hg}
There exist two positive constants $c_1,c_2$
such that, for all $x,y\in \Z^2\setminus\{0\}$,
\begin{equation}
\label{eq_order_hg}
 \frac{c_1}{\ln(1+\|x\|\vee\|y\|)} 
 \leq \hg(x,y) \leq  \frac{c_2}{\ln(1+\|x\|\vee\|y\|)}. 
\end{equation}
\end{lem}

\begin{proof}
 Assume without restricting generality that $\|x\|\geq \|y\|$
(recall that $\hg(x,y)=\hg(y,x)$),
and consider the following two cases.

\smallskip
\noindent 
Case 1: $\|y\| > \|x\|^{1/2}$. In this case~$a(x)$ and~$a(y)$
are of the same order, and, since $\|x-y\|\leq 2\|x\|$,
due to~\eqref{formula_for_a},
$a(x-y)-a(x)$ is bounded above by a positive constant;
therefore, the expression $a(x)+a(y)-a(x-y)$
will be of order~$\ln \|x\|$. This implies that~$\hg(x,y)$ will be 
of order $\frac{1}{\ln \|x\|}$ indeed.

\smallskip
\noindent 
Case 2: $\|y\| \leq \|x\|^{1/2}$. Here, \eqref{diff_a_logs}
implies that $a(x)-a(x-y) = O(\frac{\|x\|^{1/2}}{\|x\|})=O(\|x\|^{-1/2})$,
so
\[
 \hg(x,y) = \frac{a(y)+O(\|x\|^{-1/2})}{a(x)a(y)}
  = \frac{1}{a(x)} 
  \big(1+O\big(\tfrac{1}{\|x\|^{1/2}\ln(1+\|y\|)}\big)\big),
\]
and this again implies~\eqref{eq_order_hg}.
\end{proof}

It will be important to have difference estimates for 
the function~$\hg$ as well: 
\begin{lem}
\label{l_ocenka_diff_hatg}
Assume that $x,y,z\in\Z^2\setminus\{0\}$ are distinct 
and such that $\|x-y\|\wedge \|x-z\|\geq 5\|y-z\|$.
Then
\begin{equation}
\label{eq_ocenka_diff_hatg}
 \big|\hg(x,y)-\hg(x,z)\big| 
\leq O\big(\tfrac{\|y-z\|}{\|x-y\|\ln(1+\|x\|\vee\|y\|\vee\|z\|)
 \ln(1+\|y\|\vee\|z\|)}\big).
\end{equation}
\end{lem}

\begin{proof}
 First, let us write
\begin{align}
 \lefteqn{\hg(x,y)-\hg(x,z)}\nonumber\\
& = \frac{a(x)+a(y)-a(x-y)}{a(x)a(y)} 
    - \frac{a(x)+a(z)-a(x-z)}{a(x)a(z)}\nonumber\\
&= \frac{a(x)a(z)-a(x-y)a(z)-a(x)a(y)+a(x-z)a(y)}{a(x)a(y)a(z)}
   \nonumber\\
 \intertext{\footnotesize{
(put $\pm a(x-z)a(z)$ to the numerator, then
group accordingly)}}
&= \frac{a(x)(a(z)-a(y))-a(x-z)(a(z)-a(y))
 + a(z)(a(x-z)-a(x-y))}{a(x)a(y)a(z)}.
\label{oc_diff_long}
\end{align}
Throughout this proof, let us assume without loss of generality that 
$\|y\|\geq \|z\|$.
 Since the walk~$\s$ is not spatially homogeneous
(and, therefore, $\hg$ is not translationally
invariant), we need
 to take into account the relative positions of the three sites
with respect to the origin.
Specifically, we will consider the following three different cases
(see Figure~\ref{f_42threecases}).
\begin{figure}
\begin{center}
\includegraphics[width=\textwidth]{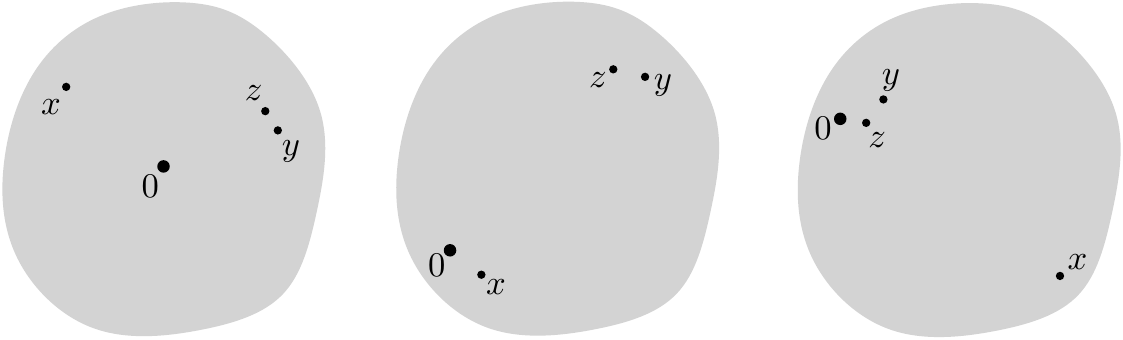}
\caption{On the proof of Lemma~\ref{l_ocenka_diff_hatg},
the three cases to consider (from left to right):
(1) $\|x\|,\|y\|$ are of the same 
\emph{logarithmic} order and $\|x\|$ is not much 
larger than $\|y\|$,
(2) $\|x\|$ is much smaller than $\|y\|$,
(3) $\|x\|$ is \emph{significantly}
larger than $\|y\|$.}
\label{f_42threecases}
\end{center}
\end{figure}

\medskip
 Case 1: $\|y\|^{1/2} \leq \|x\| \leq 2 \|y\|$.

\smallskip

In this case, the first thing to note is that 
\[
 \|y-z\| \leq \frac{\|x-y\|}{5} \leq \frac{\|x\|+\|y\|}{5}
 \leq \frac{2\|y\|+\|y\|}{5} = \frac{3}{5}\|y\|,
\]
so $\|z\|\geq \frac{2}{5}\|y\|$, meaning that~$\|y\|$
and~$\|z\|$ must be of the same order; this then implies
that $a(x),a(y),a(z)$ are all of the same order too.
Then, we use~\eqref{diff_a_logs} on the three
parentheses in the numerator of~\eqref{oc_diff_long},
to obtain after some elementary calculations
that the expression there is at most 
of order~$\frac{\|y-z\|}{\|x-y\|}\ln \|y\|$,
while the denominator is of order $\ln^3 \|y\|$.
This proves~\eqref{eq_ocenka_diff_hatg} in case~1.

\medskip
Case 2: $\|x\| < \|y\|^{1/2}$.

\smallskip
Here, it is again easy to see that~$\|y\|$
and~$\|z\|$ must be of the same order.
Now, we note that, by~\eqref{diff_a_logs},
$a(x-z)=a(z)+O\big(\tfrac{\|x\|}{\|y\|}\big)$, so, inserting 
this to~\eqref{oc_diff_long}
(and also using that $a(y)-a(z)=O\big(\tfrac{\|y-z\|}{\|y\|}\big)$),
we find that it is equal to
\begin{align}
\lefteqn{}&
 \frac{a(x)(a(z)-a(y))
 + a(z)(a(y)-a(z)+a(x-z)-a(x-y))
 +O\big(\tfrac{\|x\|}{\|y\|}\cdot\tfrac{\|y-z\|}{\|y\|}\big)}
 {a(x)a(y)a(z)}
  \nonumber\\
& = 
   \frac{a(z)-a(y)}{a(y)a(z)} 
   + \frac{a(y)-a(z)+a(x-z)-a(x-y)}{a(x)a(y)}
    + O\big(\tfrac{\|x\|\cdot \|y-z\|}{\|y\|^2\ln^2 \|y\|}\big).
\end{align}
Now, by~\eqref{diff_a_logs} the first term 
is~$O\big(\tfrac{\|y-z\|}{\|y\|\ln^2 \|y\|}\big)$ 
(that is, exactly what we need, since $\|y\|$ and $\|y-x\|$
are of the same order),
and the third term is clearly of smaller order. 
As for the second term, note that, by~\eqref{formula_for_a}
and using the fact that 
$\big|\|x-y\|\cdot\|z\|-\|y\|\cdot \|x-z\|\big|
  \leq \|x\|\cdot \|y-z\|$ by the Ptolemy's inequality,
  we obtain
\begin{align*}
\lefteqn{ a(y)-a(z)+a(x-z)-a(x-y)}\\
&= 
 \frac{2}{\pi}
 \ln \frac{\|y\|\cdot \|x-z\|}{\|x-y\|\cdot\|z\|} + O(\|y\|^{-2})\\
& =  \frac{2}{\pi} 
 \ln\Big(1 - \frac{\|x-y\|\cdot\|z\|-\|y\|\cdot \|x-z\|}
 {\|x-y\|\cdot\|z\|} \Big)+ O(\|y\|^{-2})\\
 & =  O\big(\tfrac{\|x\|\cdot \|y-z\|}{\|x-y\|\cdot\|z\|}
 +\|y\|^{-2}\big) = 
  O\big(\tfrac{\|x\|\cdot \|y-z\|}{\|z\|^2}\big),
\end{align*}
so it is again of smaller order than the first term.
This shows~\eqref{eq_ocenka_diff_hatg} in case~2.

\medskip
Case 3: $\|x\| > 2 \|y\|$.

\smallskip

Notice that, in this case, $\|z\|$ need not be of the 
same order as~$\|y\|$, it may happen to be significantly smaller.
Here (by also grouping the first two terms in the numerator)
we rewrite~\eqref{oc_diff_long} as
\begin{equation}
\label{rewrite_case3}
\frac{(a(x)-a(x-z))(a(z)-a(y))}{a(x)a(y)a(z)}
  + \frac{a(x-z)-a(x-y)}{a(x)a(y)}.
\end{equation}
By~\eqref{diff_a_logs}, the second term 
is~$O\big(\tfrac{\|y-z\|}{\|x-y\|\ln (1+\|x\|)\ln(1+\|y\|)}\big)$ 
(that is, exactly what we need). 
Next, observe that (recall that we assumed that $\|y\|\geq\|z\|$)
\[
\ln \|y\| - \ln \|z\| =
\ln\frac{\|y\|}{\|z\|} \leq \frac{\|y\|}{\|z\|} -1
 = \frac{\|y\|-\|z\|}{\|z\|} \leq  \frac{\|y-z\|}{\|z\|}.
\]
Therefore (also using~\eqref{diff_a_logs} on the first factor),
the numerator of the first term is 
$O\big(\tfrac{\|z\|}{\|x\|}\times\tfrac{\|y-z\|}{\|z\|}\big)
 = O\big(\tfrac{\|y-z\|}{\|x\|}\big)
 $, and 
 so (since the denominator is not less than~$a(x)a(y)$) 
 the first term in~\eqref{rewrite_case3}
is at most of the same order as the second one. 
This concludes the proof of Lemma~\ref{l_ocenka_diff_hatg}.
\end{proof}

Lemma~\ref{l_ocenka_diff_hatg} permits us to obtain 
the following useful
 expression for the probability of ever hitting~$A$
 from a distant site.
Consider a finite~$A\subset \Z^2\setminus\{0\}$
with $y_0\in A$ and note first that $\|y_0\|+\diam(A)$
will be of the same order regardless of the choice of~$y_0$. 
Assume now that $\dist(x,A)>5\diam(A)$;
then, it also holds that $\|x\|\vee\|y\|\vee\|z\|$
is of the same order as $\|x\|\vee (\|y_0\|+\diam(A))$
for any choice of $y,z\in A$. Indeed, trivially
$\|x\|\vee\|y\|\vee\|z\|\leq \|x\|\vee (\|y_0\|+\diam(A))$;
on the other hand, $\|y\|+\diam(A)<\|y\|+\dist(x,A)
\leq \|y\|+\|x-y\|\leq \|x\|+2\|y\| \leq 3(\|x\|\vee\|y\|\vee\|z\|)$.
Then, note that $1+\|y\|+\|y-z\|\leq 1+2\|y\|+\|z\| \leq 3(1+\|y\|\vee\|z\|)
<(1+\|y\|\vee\|z\|)^3$ (since the expression in parentheses is at least~$2$),
so 
\[
 \frac{\|y-z\|}{\ln (1+\|y\|\vee\|z\|)} \leq 
  3 \frac{\|y-z\|}{\ln (1+\|y\| +\|y-z\|)} 
   \leq 3 \frac{\diam(A)}{\ln (1+\|y\| + \diam(A))},
\]
where the second inequality is
due to the fact that the function $f(x)=\frac{x}{\ln(a+x)}$
is increasing on $(0,+\infty)$ for any $a\geq 1$.
So, recalling~\eqref{df_hcapa} and~\eqref{Green_hEs},
 we see that Lemma~\ref{l_ocenka_diff_hatg}
 gives us (in the case $\dist(x,A)>5\diam(A)$) that
\begin{equation}
\label{form_nothitA_hS}
 \IP_x[\htau_A<\infty] =\hcapa(A)\big(\hg(x,y_0)
+ O\big(\tfrac{\diam(A)}{\dist(x,A)
\ln(1+\|x\|\vee (\|y_0\|+\diam(A)))
\ln(1+(\|y_0\|+\diam(A)))}\big)\big).
\end{equation}

The next technical fact that we need is that
the conditioned walk can go out of an annulus
with uniformly positive probability:
%
\begin{lem}
\label{l_diff_hS}
 Let $b,C$ be positive constants such that
$1+b<C$, and assume that $r\geq 1, x_0,y_0\in\Z^2\setminus\{0\}$
are such that $x_0 \in \B(y_0, Cr)$ and $\|x_0-y_0\| > (1+b)r$.
Then, there exists a constant~$c'>0$ (depending \emph{only}
on~$b$ and~$C$)
such that
\begin{equation}
\label{eq_diff_hS}
 \IP_{x_0}[\htau_{\partial\B(y_0,Cr)}<\htau_{\B(y_0,r)}]\geq c'.
\end{equation}
\end{lem}

\begin{proof}
Note that we can assume that~$r$ is large enough, otherwise
the uniform ellipticity of the conditioned walk 
will imply the claim.
First of all, it is clear 
that the analogue of~\eqref{eq_diff_hS}
holds for simple random walk, i.e.,
for all $r,x_0,y_0$ as above it holds that 
\begin{equation}
\label{eq_diff_SRW}
 \IP_{x_0}[\tau_{\partial\B(y_0,2Cr)}<\tau_{\B(y_0,r)}]\geq c''
\end{equation}
for some $c''>0$ which depends only on~$b$ and~$C$.
Now, the idea is to derive~\eqref{eq_diff_hS}
 from~\eqref{eq_diff_SRW}. It holds that
 the weight~$\hP_\vr$ 
with respect to~$\s$
 of a finite path~$\vr$ which starts at~$x_0$ and
does not pass through the origin equals 
$\frac{a(\vr_{\text{end}})}{a(x_0)}P_\vr$, where 
$P_\vr=(\frac{1}{4})^{|\vr|}$ (with~$|\vr|$ being the length
of the path~$\vr$ and $\vr_{\text{end}}$ being the last site of it)
is the weight of the same path with respect to the simple random walk.
We can then write for any~$R>0$ such that $x_0\in\B(y_0,R)$
\begin{align}
\lefteqn{
\IP_{x_0}[\htau_{\partial\B(y_0,R)}<\htau_{\B(y_0,r)}]
}\nonumber\\
    &\geq 
 \min_{z\in\partial\B(y_0,R)}\frac{a(z)}{a(x_0)}\times
 \IP_{x_0}[\tau_{\partial\B(y_0,R)}<\tau_{\B(y_0,r)},
   \tau_{\partial\B(y_0,R)}<\tau_0]\nonumber\\
 & \geq 
  \min_{z\in\partial\B(y_0,R)}\frac{a(z)}{a(x_0)}\times
 \Big( \IP_{x_0}[\tau_{\partial\B(y_0,R)}<\tau_{\B(y_0,r)}]
    - \IP_{x_0}[\tau_0< \tau_{\partial\B(y_0,R)}]\Big).
\label{main_oc_diff_hS}    
\end{align}
\begin{figure}
\begin{center}
\includegraphics{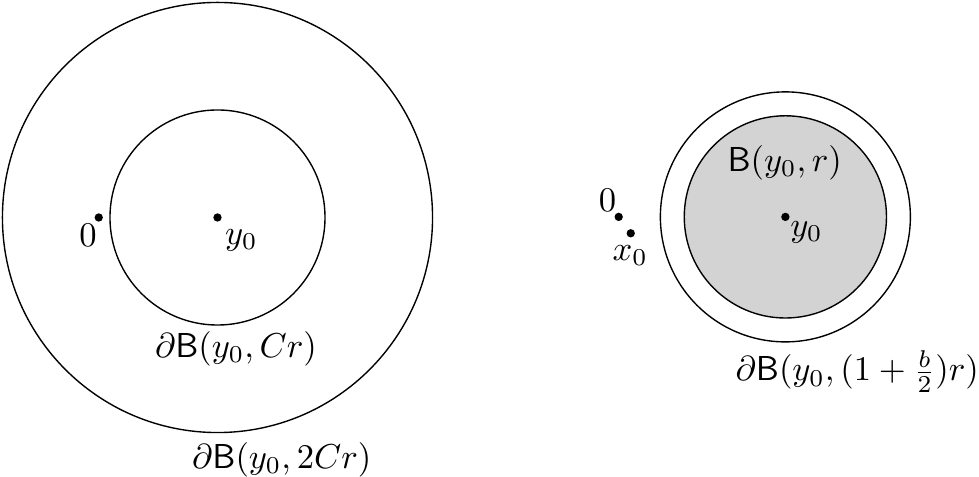}
\caption{On the proof of Lemma~\ref{l_diff_hS}}
\label{f_43diff_hS}
\end{center}
\end{figure}
 Now, a key observation is that, for all large enough~$r$,
the property
\begin{equation}
\label{assump_on_B(y0,R)} 
 \frac{a(z)}{a(x_0)}\geq \frac{1}{2}
\text{ for all } z\in\partial \B(y_0,R)
\end{equation}
holds for either $R=Cr$ or $R=2Cr$ (or both). Indeed,
roughly speaking, for~\eqref{assump_on_B(y0,R)} to hold 
it would be enough that~$\|z\|$ is of order~$r+\|y_0\|$
for all $z\in\partial \B(y_0,R)$; this can be seen to be so
in at least one of the above cases
 (look at the left side
of Figure~\ref{f_43diff_hS}: if~$\partial\B(y_0,Cr)$ is ``too 
close'' to the origin, then~$\partial\B(y_0,2Cr)$ is not).

For definiteness, assume now that~\eqref{assump_on_B(y0,R)} 
holds for $R=2Cr$. 
By~\eqref{main_oc_diff_hS}, we have then
\begin{align*}
\lefteqn{
\IP_{x_0}[\htau_{\partial\B(y_0,Cr)}<\htau_{\B(y_0,r)}]
}\\
 & \geq \IP_{x_0}[\htau_{\partial\B(y_0,2Cr)}<\htau_{\B(y_0,r)}]\\
 & \geq 
  \frac{1}{2}\Big( \IP_{x_0}[\tau_{\partial\B(y_0,2Cr)}<\tau_{\B(y_0,r)}]
    - \IP_{x_0}[\tau_0< \tau_{\partial\B(y_0,2Cr)}]\Big)\\
 &\geq  \frac{1}{2}\Big( c''
 - \IP_{x_0}[\tau_0< \tau_{\partial\B(y_0,2Cr)}]\Big)\\
 &\geq \frac{c''}{4},
\end{align*}
provided that 
\begin{equation}
\label{P_to0_small}
 \IP_{x_0}[\tau_0< \tau_{\partial\B(y_0,2Cr)}] \leq \frac{c''}{2}.
\end{equation}
Now, if $0\notin \B(y_0,Cr)$, then~\eqref{P_to0_small}
trivially holds; so, let us assume that $0\in \B(y_0,Cr)$.
We then consider two cases: $\|x_0\|\geq \frac{b}{4}r$,
and $\|x_0\| < \frac{b}{4}r$.
In the first case, Lemma~\ref{l_escape_origin_d=2}
implies that
$\IP_{x_0}[\tau_0< \tau_{\partial\B(y_0,2Cr)}]\asymp \frac{1}{\ln r}$,
so~\eqref{P_to0_small} holds for large enough~$r$.
In the second case, note first that 
\[
 \IP_{x_0}[\htau_{\partial\B(y_0,Cr)}<\htau_{\B(y_0,r)}]
   \geq \min_{z\in \partial\B(y_0, (1+\frac{b}{2})r)}
     \IP_z[\htau_{\partial\B(y_0,Cr)}<\htau_{\B(y_0,r)}]
\]
and, for all $z\in \partial\B(y_0, (1+\frac{b}{2})r)$,
it holds that $\|z\|\geq \frac{b}{4}r$
(see Figure~\ref{f_43diff_hS} on the right).
We may then repeat the above reasoning with an arbitrary 
$z\in \partial\B(y_0, (1+\frac{b}{2})r)$ on the place 
of~$x_0$ to finally obtain the claim.
\end{proof}

Next technical fact we need is the following lower
bound on the probability that the conditioned walk
never hits a disk:
\begin{lem}
\label{l_escapeanydisk_hS}
Fix~$b>0$ and assume 
 that $x_0,y_0\in\Z^2\setminus\{0\}$ and $r\geq 1$
are such that $\|x_0-y_0\|\geq (1+b) r$.
Then there exists $c=c(b)$ such that for all $x_0,y_0,r$ as above
it holds that
\begin{equation}
\label{eq_escapeanydisk_hS} 
\IP_{x_0}[\htau_{\B(y_0,r)}=\infty] 
  \geq \frac{c}{\ln (\|y_0\|+r)}.
\end{equation}
\end{lem}
\begin{proof}
 We need to consider two cases: $\B(y_0,r)$
is (\emph{relatively} to its size) close to/far from the origin.
First, let us assume that $\|y_0\|< 12r$ (so that the disk
is relatively close to the origin).
In this case, we can assume additionally that $\|x_0-y_0\|\geq 51 r$
(indeed, 
Lemma~\ref{l_diff_hS} implies that, for any starting position~$x_0$
such that $\|x_0-y_0\|\geq (1+b) r$,
with at least a constant
probability the walk reaches $\partial\B(y_0,51r)$ before 
hitting~$\B(y_0,r)$).
Then, it holds that $\B(y_0,r)\subset\B(13r)$,
and 
$\B(26r)\subset \B(y_0, 51r)$.
Now, Lemma~\ref{l_escape_hatS} easily implies
that, if $r\geq 1$ and $\|x\|\geq 2r$
\[
 \IP_x[\htau_{\B(r)}=\infty] \geq \frac{c'}{\ln r}
\]
(because~\eqref{eq_escape_hatS} will work for 
large enough~$\|x\|$, and one can use the uniform
ellipticity of~$\s$ otherwise); 
this proves~\eqref{eq_escapeanydisk_hS} in the first
case.

Now, suppose that $\|y_0\|\geq 12r$ (that is,
$r\leq \frac{1}{12}\|y_0\|$).
Analogously to the previous case, 
Lemma~\ref{l_diff_hS} permits us to
assume without loss of generality that $x_0\in \partial \B(y_0,3r)$.
We now use 
the martingale (recall Proposition~\ref{p_hG_mart})
\[
 \hl(\s_{n\wedge \htau_{y_0}},y_0) 
   = 1+\frac{a(y_0)-a(\s_{n\wedge \htau_{y_0}}-y_0)}
   {a(\s_{n\wedge \htau_{y_0}})}.
\]
The Optional Stopping Theorem implies that
\begin{align*}
\hl(x_0,y_0) 
&= \sum_{z\in \partial \B(y_0,r)} 
 \IP_{x_0}[\htau_{\B(y_0,r)}<\infty, \s_{\htau_{\B(y_0,r)}}=z]
  \hl(z,y_0)\\
& \geq \IP_{x_0}[\htau_{\B(y_0,r)}<\infty]
 \min_{z\in \partial \B(y_0,r)} \hl(z,y_0),
\end{align*}
so
\begin{equation}
\label{oc_ubezhal_B_y_0_r}
 \IP_{x_0}[\htau_{\B(y_0,r)}<\infty]
\leq \frac{\hl(x_0,y_0)}{\min_{z\in \partial \B(y_0,r)}\hl(z,y_0)}.
\end{equation}

Assume $z\in \partial \B(y_0,r)$ and write, 
using~\eqref{formula_for_a} and with $\gamma'':=\pi\gamma'/2$,
\begin{align}
 \hl(z,y_0) &= \frac{a(z)+a(y_0)-a(y_0-z)}{a(z)}\nonumber\\
 &= \frac{\ln\|z\| +\ln\|y_0\| - \ln r + \gamma''+ O(\|y_0\|^{-2}+r^{-1})}
 {\ln\|z\| + \gamma'' + O(\|z\|^{-2})} \nonumber\\
 &\geq \frac{\ln(\|y_0\|-r) + 
 \ln\|y_0\| - \ln r + \gamma''+ O(\|y_0\|^{-2}+r^{-1})}
 {\ln(\|y_0\|+r) + \gamma'' + O(\|y_0\|^{-2}+r^{-1})} \nonumber\\
 &= \frac{ 2\ln\|y_0\|
 +\ln\big(1-\frac{r}{\|y_0\|}\big) 
 - \ln r + \gamma''+ O(r^{-1})}
 {\ln\|y_0\|+\ln\big(1+\frac{r}{\|y_0\|}\big) 
 + \gamma'' + O(r^{-1})} \nonumber\\
 & := \frac{T_1}{T_2},
\label{hl_z_geq}
\end{align}
and, denoting $R:=\|x_0-y_0\|$ (so that $R=3r+O(1)$),
\begin{align}
 \hl(x_0,y_0) &= \frac{a(x_0)+a(y_0)-a(y_0-x_0)}{a(x_0)}\nonumber\\
 &= \frac{\ln\|x_0\| +\ln\|y_0\| - \ln R + \gamma''+ O(\|y_0\|^{-2}+R^{-2})}
 {\ln\|x_0\| + \gamma'' + O(\|x_0\|^{-2})} \nonumber\\
 &\leq \frac{\ln(\|y_0\|+R)
 + \ln\|y_0\| - \ln R + \gamma''+ O(\|y_0\|^{-2}+R^{-2})}
 {\ln(\|y_0\|-R) + \gamma'' + O(\|y_0\|^{-2}+R^{-2})} \nonumber\\
& =  \frac{ 2\ln\|y_0\|
 +\ln\big(1+\frac{R}{\|y_0\|}\big) 
 - \ln R + \gamma''+ O(R^{-2})}
 {\ln\|y_0\|+\ln\big(1-\frac{R}{\|y_0\|}\big) 
 + \gamma'' + O(R^{-2})} \nonumber\\
  & := \frac{T_3}{T_4} .
\label{hl_x_geq}
\end{align}
Now, a straightforward calculation yields 
\begin{align*}
 \frac{T_2}{T_4} & = 1 + \frac{\ln\frac{1+r/\|y_0\|}{1-R/\|y_0\|}
  + O(r^{-1})}{\ln\|y_0\|+\ln\big(1-\frac{R}{\|y_0\|}\big) 
 + \gamma'' + O(R^{-2})},\\
 \intertext{and}
  \frac{T_3}{T_1} & = 
 1 - \frac{\ln\big(\frac{R}{r}\cdot 
  \frac{1-r/\|y_0\|}{1+R/\|y_0\|}\big)
  + O(r^{-1})}{2\ln\|y_0\|
  -\ln\frac{r}{1-r/\|y_0\|} + \gamma''+ O(r^{-1})}\\
  &\leq  1 - \frac{\ln\big(\frac{R}{r}\cdot 
  \frac{1-r/\|y_0\|}{1+R/\|y_0\|}\big)
  + O(r^{-1})}{2\ln\|y_0\|+ \gamma''}.
\end{align*}
Therefore, 
by~\eqref{oc_ubezhal_B_y_0_r} 
we have (after some elementary calculations)
\begin{align}
\lefteqn{
\IP_{x_0}[\htau_{\B(y_0,r)}<\infty]} \nonumber\\
   & \leq \frac{T_2}{T_4} \times \frac{T_3}{T_1} \nonumber\\
   & \leq 1 - \frac{\ln\big(\frac{R}{r}\cdot 
  \frac{1-r/\|y_0\|}{1+R/\|y_0\|}\big)^{1/2}
  - \ln \frac{1+r/\|y_0\|}{1-R/\|y_0\|}+O((\ln r)^{-1})}
  {\ln\|y_0\|
  \big(1+ O\big(\tfrac{1}{\ln\|y_0\|}\big)\big)}.
\label{1-const/ln}
\end{align}
It remains only to observe that, if~$r$
is large enough, the numerator in~\eqref{1-const/ln}
is bounded from below by a positive constant: indeed,
observe that $\frac{R}{r}$ is (asymptotically)
 $3$, $\frac{r}{\|y_0\|}$ and 
$\frac{R}{\|y_0\|}$ are at most~$\frac{1}{12}$
and~$\frac{1}{4}$ respectively, and 
\[
 \sqrt{3 \times \frac{1-\frac{1}{12}}{1+\frac{1}{4}}}
  \times \frac{1-\frac{1}{4}}{1+\frac{1}{12}}
  =\sqrt{\frac{891}{845}} > 1.
\]
This concludes the proof of Lemma~\ref{l_escapeanydisk_hS}
in the case when~$r$ is large enough; the case of smaller
values of~$r$, though, can be easily reduced
to the former one by using the uniform ellipticity
of the $\s$-walk.
\end{proof}

\section{Proofs of the main results}
\label{s_proofs}

\begin{proof}[Proof of Theorem~\ref{t_GF_conditional}]
First, we need a very simple general fact about hitting times 
of recurrent Markov chains:

\begin{lem}
\label{l_escape_x_AB}
Let $(X_n)$ be a recurrent Markov chain on a state space~$\Sigma$,
and $x\in \Sigma$, $A,B \subset \Sigma$ are such that
$A\cap B = \emptyset$ and $x\notin A\cup B$. Then
\begin{equation}
\label{eq_escape_x_AB} 
 \IP_x[\tau_A<\tau_B] 
   = \IP_x[\tau_A<\tau_B\mid \tau^+_x > \tau_{A\cup B}] 
\end{equation}
(that is, the events $\{\tau_A<\tau_B\}$ and $\{\tau^+_x > \tau_{A\cup B}\}$
are independent under~$\IP_x$).
\end{lem}

\begin{proof}
This is almost evident, so we give only a sketch
(see Figure~\ref{f_escape_x_AB}): 
let $p:=\IP_x[\tau_A<\tau_B\mid \tau^+_x > \tau_{A\cup B}]$
be the value of the probability in the right-hand side 
of~\eqref{eq_escape_x_AB}.
At the moments when the walker visits~$x$, it tosses a coin to 
decide if it will revisit it before coming to~$A\cup B$, or not.
When it decides to definitely leave~$x$ for~$A\cup B$, the 
probability of choosing~$A$ is~$p$, so it is~$p$ overall.
\begin{figure}
\begin{center}
\includegraphics{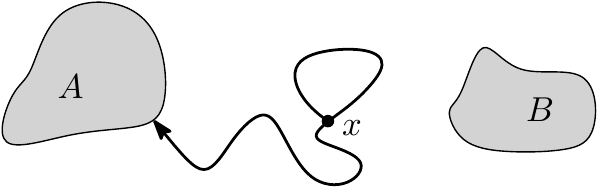}
\caption{On the proof of Lemma~\ref{l_escape_x_AB}}
\label{f_escape_x_AB}
\end{center}
\end{figure}
\end{proof}
We continue proving Theorem~\ref{t_GF_conditional}.
Fix (a large) $R>0$, abbreviate~$\Lambda_R=\B(R)\setminus\{0\}$, 
and let us denote for $y\in \Lambda_R$
\begin{align*}
N^*_{y,R} &= \sum_{k=0}^{\tau_{\Lambda_R^c}} \1{S_k=y},\\
\hN^*_{y,R} &= \sum_{k=0}^{\htau_{\Lambda_R^c}} \1{\s_k=y},
\end{align*}
to be the numbers of visits to~$y$ before hitting~$0$
or~$\partial_e\B(R)$,
for the simple random walk and the conditioned walk.
Let us also denote $\hG_R(x,y)=\E_x \hN^*_{y,R}$.
Now, let $x\in \Lambda_R$ and observe that, on one hand,
\begin{align}
\lefteqn{
\IP_x[N^*_{y,R}=n,\tau_{\partial_e\B(R)}<\tau_0]
}\nonumber\\
 &= \IP_x[N^*_{y,R}\geq n]
  \IP_y[\tau_{\partial_e\B(R)}<\tau_0,\tau^+_y> \tau_{\Lambda_R^c}]
  \nonumber\\
 &= \IP_x[N^*_{y,R}\geq n]
  \IP_y[\tau^+_y> \tau_{\Lambda_R^c}]
  \IP_y[\tau_{\partial_e\B(R)}<\tau_0 \mid \tau^+_y> \tau_{\Lambda_R^c}] 
  \nonumber\\
  \intertext{\footnotesize{\qquad \qquad (by Lemma~\ref{l_escape_x_AB}) }}
 &= \IP_x[N^*_{y,R}\geq n]
  \IP_y[\tau^+_y> \tau_{\Lambda_R^c}]
  \IP_y[\tau_{\partial_e\B(R)}<\tau_0 ] 
  \nonumber\\ 
 &= \IP_x[N^*_{y,R}=n]
  \IP_y[\tau_{\partial_e\B(R)}<\tau_0 ]  \nonumber\\ 
     \intertext{\footnotesize{\qquad \qquad 
     (by Lemma~\ref{l_escape_origin_d=2}) }}
 &= \IP_x[N^*_{y,R}=n]
  \frac{a(y)}{a(R)+O(R^{-1})},
 \label{in_form_G}
\end{align}
and, on the other hand, the same expression 
can be also treated in the following way:
\begin{align}
\lefteqn{
\IP_x[N^*_{y,R}=n,\tau_{\partial_e\B(R)}<\tau_0]
}\nonumber\\
 &= \IP_x[N^*_{y,R}=n\mid\tau_{\partial_e\B(R)}<\tau_0]
   \IP_x[\tau_{\partial_e\B(R)}<\tau_0]\nonumber\\
       \intertext{\footnotesize{\qquad \qquad
       (by Lemma~\ref{l_relation_S_hatS}) }}
 &= \IP_x[\hN^*_{y,R}=n] \big(1+O((R\ln R)^{-1})\big) 
 \IP_x[\tau_{\partial_e\B(R)}<\tau_0]   \nonumber \\
      \intertext{\footnotesize{\qquad \qquad
      (again, by Lemma~\ref{l_escape_origin_d=2}) }}
 &= \IP_x[\hN^*_{y,R}=n]\big(1+O((R\ln R)^{-1})\big) 
  \frac{a(x)}{a(R)+O(R^{-1})}.
   \label{in_form_hG}
\end{align}
Note also that $a(R)+O(R^{-1})=a(R)\big(1+O((R\ln R)^{-1})\big)$.
So, since~\eqref{in_form_G} and~\eqref{in_form_hG}
must be equal, we have
\[
a(x)\IP_x[\hN^*_{y,R}=n]
  = a(y)\IP_x[N^*_{y,R}=n]\big(1+O((R\ln R)^{-1})\big) ;
\]
multiplying by~$n$ and summing in~$n\geq 1$, we obtain
\begin{equation}
\label{G_hG} 
 a(x)\hG_R(x,y) = a(y)G_{\Lambda_R}(x,y)\big(1+O((R\ln R)^{-1})\big).
\end{equation}
Note that $\hG_R(x,y)\to \hG(x,y)$ as $R\to\infty$,
due to the Monotone Convergence Theorem.
Next, we are actually able to say something
about~$G_{\Lambda_R}(x,y)$:
by Proposition~4.6.2(b) of~\cite{LawLim10}
it holds that
\begin{align*}
 G_{\Lambda_R}(x,y) &= \IE_x a(S_{\tau_{\Lambda_R^c}}-y)-a(x-y)\\
       \intertext{\footnotesize{\qquad\qquad \qquad
      (once again, by Lemma~\ref{l_escape_origin_d=2}) }}
 &= \frac{a(x)}{a(R)+O(R^{-1})} \big(a(R)
  +O\big(\tfrac{\|y\|+1}{R}\big)\big)\\
 &\qquad
   + \Big(1-\frac{a(x)}{a(R)+O(R^{-1})}\Big)a(y)-a(x-y)\\
 &= a(x) + a(y)-a(x-y) 
 + O\big(\tfrac{\|y\|+1}{R} + \tfrac{a(x)a(y)}{a(R)}\big).
\end{align*}
Inserting this back to~\eqref{G_hG} and sending~$R$
to infinity, we finally obtain~\eqref{eq_GF_conditional}.
\end{proof}

\begin{proof}[Proof of Theorem~\ref{t_conv_harm_hatS}]
Let us assume without restricting generality 
that~$A$ contains at least two sites,
so that $\diam(A)\geq 1$.
For $x\notin A$, $y\in \partial A$, and $n\geq 1$,
let us denote by~$\Theta_{xy}^{(n)}$  the set of nearest-neighbour
trajectories $\wp=(z_0,\ldots,z_k)$ such that
\begin{itemize}
 \item $z_0=x$, $z_k=y$, and $z_j\notin A$ for all $j\leq k-1$,
 i.e., the trajectory ends on the first entrance to~$A$, which takes
place in~$y$;
 \item $\sum_{j=0}^{k}\1{z_j=x}=n$, i.e., the trajectory
 visits~$x$ exactly~$n$ times (note that we \emph{do} count
$z_0=x$ as one visit);
\end{itemize}
 Let us also denote by
\[
 \hN_x = \sum_{j=0}^\infty \1{\s_j=x}
\]
the total number of visits to~$x\notin A$,
by
\[
 \hN_x^{\flat} = \sum_{j=0}^{\htau^+_A-1} \1{\s_j=x}
\]
the number of visits to~$x$ \emph{before} the first return to~$A$,
and by
\[
\hN_x^{\sharp} = \sum_{j=\htau^+_A}^\infty \1{\s_j=x}
\]
the number of visits to~$x$ \emph{after} the first return to~$A$
(naturally, setting $\hN_x^{\sharp}=0$ on $\{\htau^+_A=\infty\}$).

Recall that
$P_\wp$ (respectively, $\hP_\wp$)
is the weight of the trajectory~$\wp$
with respect to the simple random walk (respectively, to the $\s$-walk).
First,
it is clear that
\[
  \IP_x[ \htau_A<\infty, \s_{\htau_A}=y] 
    = \sum_{n=1}^{\infty} \sum_{\wp\in \Theta_{xy}^{(n)}}\hP_{\wp};
\]
then, by the reversibility of the $\s$-walk,
\[
 \IP_y[\hN_x^{\flat}\geq n] 
  \frac{a^2(y)}{a^2(x)}= \sum_{\wp\in \Theta_{xy}^{(n)}}\hP_{\wp}.
\]
Now, we can write
\begin{align}
\lefteqn{ \IP_x[\s_{\htau_A}=y\mid \htau_A<\infty] }\nonumber\\
 &= \frac{\IP_x[ \htau_A<\infty, \s_{\htau_A}=y]}{\IP_x[ \htau_A<\infty]}
 \nonumber\\
 &= \frac{1}{\IP_x[ \htau_A<\infty]}\sum_{n=1}^{\infty} 
 \frac{a^2(y)}{a^2(x)}\IP_y[\hN_x^{\flat}\geq n]
 \nonumber\\
 &= \frac{a^2(y)}{a^2(x)\IP_x[ \htau_A<\infty]}
 \IE_y \hN_x^{\flat}
 \nonumber\\
&=  \frac{a^2(y)}{a^2(x)\IP_x[ \htau_A<\infty]}
(\IE_y \hN_x - \IE_y \hN_x^{\sharp})
\nonumber\\
&= \frac{a^2(y)}{\IP_x[ \htau_A<\infty]}
\Big(\frac{\hG(y,x)}{a^2(x)}
 - \sum_{z\in \partial A}
  \IP_y[ \htau^+_A<\infty, \s_{\htau^+_A}=z]\frac{\hG(z,x)}{a^2(x)}\Big)
 \nonumber\\
 &= \frac{a^2(y)}{\IP_x[ \htau_A<\infty]}
 \Big(\hg(y,x)
 - \sum_{z\in \partial A}
  \IP_y[ \htau^+_A<\infty, \s_{\htau^+_A}=z]\hg(z,x)\Big) \nonumber\\
  &=  \frac{a^2(y)}{\IP_x[ \htau_A<\infty]}
 \Big(\hg(y,x)\Big(\hEs_A(y)+ \sum_{z\in \partial A}
  \IP_y[ \htau^+_A<\infty, \s_{\htau^+_A}=z]\Big)\nonumber\\
 &\qquad \qquad \qquad \qquad\qquad - \sum_{z\in \partial A}
  \IP_y[ \htau^+_A<\infty, \s_{\htau^+_A}=z]\hg(z,x)\Big)\nonumber\\
  &= \frac{a^2(y)\hg(y,x)\hEs_A(y)}{\IP_x[ \htau_A<\infty]}\nonumber\\
 &\quad + \frac{a^2(y)}{\IP_x[ \htau_A<\infty]}
  \sum_{z\in \partial A}
  \IP_y[ \htau^+_A<\infty, \s_{\htau^+_A}=z](\hg(y,x)-\hg(z,x)).
 \label{sdfgh345}
\end{align}

By~\eqref{form_nothitA_hS} and Lemma~\ref{l_order_hg} it holds that
\begin{equation}
\label{1st_term_bigcalc}
 \frac{a^2(y)\hg(y,x)\hEs_A(y)}{\IP_x[ \htau_A<\infty]}
  = \hhm_A(y)
\big(1+ O\big(\tfrac{\diam(A)}{\dist(x,A)}\big)\big).
\end{equation}
Therefore, it only remains to show that the second term 
in~\eqref{sdfgh345} is 
$O\big(\hhm_A(y)\tfrac{\diam(A)}{\dist(x,A)}\big)$.

%
Now,
we are going to use the Optional Stopping Theorem
with the martingale~$\hM_{n\wedge \htau_x}$ where
\[
\hM_n = \hg(y,x) - \hg(\s_n,x)
\]
to estimate the second term in~\eqref{sdfgh345}.
Recall that $y\in\partial A$, and let us define
\[
 V = \partial\B(y, 2\diam(A)),
\]
see Figure~\ref{f_32AVx_no_0}.
\begin{figure}
\begin{center}
\includegraphics{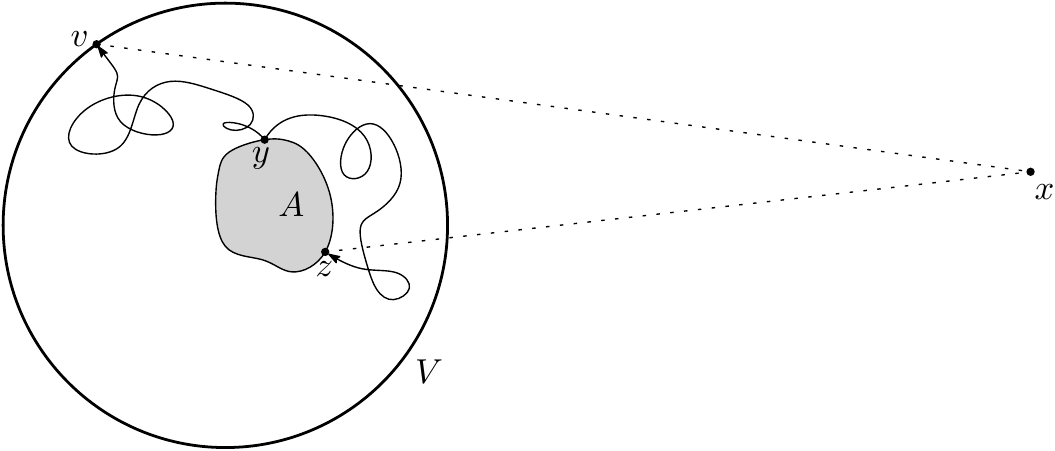}
\caption{On the proof of Theorem~\ref{t_conv_harm_hatS}:
the random walk starts at~$y\in\partial A$ and ends either
on the first re-entry to~$A$ or entry to~$V$.
}
\label{f_32AVx_no_0}
\end{center}
\end{figure}
Let $\tau = \htau^+_A\wedge \htau_V$.
We have (note that $\tau<\htau_x$)
\begin{align*}
 0 &= \IE_y \hM_0 \\
 &= \IE_y \hM_\tau \\
 &= \IE_y \big( \hM_{\htau^+_A}\1{\htau^+_A<\htau_V}\big)
  + \IE_y \big( \hM_{\htau_V}\1{\htau_V<\htau^+_A}\big)\\
    \intertext{\footnotesize{ \quad
(since $\1{\htau^+_A<\infty} = \1{\htau^+_A<\htau_V}
+ \1{\htau_V<\htau^+_A<\infty}$)}}
& = \IE_y \big( \hM_{\htau^+_A}\1{\htau^+_A<\infty}\big)
- \IE_y \big( \hM_{\htau^+_A}\1{\htau_V<\htau^+_A<\infty}\big)\\
& \qquad+\IE_y \big( \hM_{\htau_V}\1{\htau_V<\htau^+_A}\big).
\end{align*}
Note that for any $z\in V\cup \partial A$ it holds
that $\|y-z\|\leq 2\diam(A)$ and $\|x-z\|\geq 10(\diam(A)+1)$,
so in the following we will be able to apply 
Lemma~\ref{l_ocenka_diff_hatg}.
Since
\[
\IE_y\big(\hM_{\htau^+_A}\1{\htau^+_A<\infty}\big)
 = \sum_{z\in \partial A}
  \IP_y[ \htau^+_A<\infty, \s_{\htau^+_A}=z]
  \big(\hg(y,x)-\hg(z,x)\big),
\]
we obtain that 
\begin{align}
\lefteqn{ \sum_{z\in \partial A}
  \IP_y[ \htau^+_A<\infty, \s_{\htau^+_A}=z](\hg(y,x)-\hg(z,x))
} \nonumber \\
 & = \IE_y \big( \hM_{\htau^+_A}\1{\htau_V<\htau^+_A<\infty}\big)
 - \IE_y \big( \hM_{\htau_V}\1{\htau_V<\htau^+_A}\big)
 \nonumber \\
 &= \IP_y[\htau_V<\htau^+_A]
  \Big(\IE_y \big( \hM_{\htau^+_A}\1{\htau^+_A<\infty}
   \mid \htau_V<\htau^+_A\big) 
   -  \IE_y \big( \hM_{\htau_V}\mid \htau_V<\htau^+_A\big)\Big)
  \nonumber \\  
\intertext{\footnotesize{(by Lemma~\ref{l_ocenka_diff_hatg}
and the derivation of~\eqref{form_nothitA_hS};
recall that $\hM_\tau = \hg(y,x) - \hg(z,x)$
on $\{\s_{\tau}=z\}$)}}  
& \leq \IP_y[\htau_V<\htau^+_A]
 \times O\big(\tfrac{\diam(A)}{\dist(x,A)
 \ln(1+\|y\|+\diam(A))\ln(1+\|x\|\vee (\|y\|+\diam(A)))}\big).
\label{usamosL4.13}
\end{align}
Next, we can write
\begin{align*}
 \hEs_A(y) & = \IP_y[\htau^+_A=\infty]\\
 &= \sum_{v\in V} \IP_y[\htau_V<\htau^+_A,
 \s_{\htau_V}=v]\IP_v[\htau_A=\infty]\\
 \intertext{\qquad \qquad \qquad \qquad
 \footnotesize{(by Lemma~\ref{l_escapeanydisk_hS})}} 
 &\geq \frac{c}{\ln(\|y\|+\diam(A))} \IP_y[\htau_V<\htau^+_A],
\end{align*}
which means that 
\begin{equation}
\label{hS_escape_to_V}
 \IP_y[\htau_V<\htau^+_A]
 \leq O\big(\hEs_A(y)\ln(\|y\|+\diam(A)) \big).
\end{equation}
Also, \eqref{1st_term_bigcalc} implies that
\begin{equation}
\label{2.1st_term_bigcalc}
 \frac{a^2(y)}{\IP_x[ \htau_A<\infty]}
  = O\Big(\frac{\hhm_A(y)}{\hEs_A(y)\hg(x,y)}\Big).
\end{equation}
Since, by Lemma~\ref{l_order_hg},  
\[
\frac{1}{\hg(x,y)}=O\big(\ln(1 + \|x\|\vee\|y\|)\big),
\]
it only remains to combine~\eqref{usamosL4.13}%
--\eqref{2.1st_term_bigcalc}
to see that the second term 
in~\eqref{sdfgh345} is indeed 
$O\big(\hhm_A(y)\tfrac{\diam(A)}{\dist(x,A)}\big)$,
thus concluding the proof of Theorem~\ref{t_conv_harm_hatS}.
\end{proof}

\section*{Acknowledgments}
This work was partially supported by CNPq
(301605/2015--7).
The author is grateful to the anonymous referee 
for many comments and suggestions on the first version
of this paper.

\end{document}